\documentclass [12pt,a4paper,reqno]{amsart}
\usepackage{amsmath,amsthm,amscd}
\usepackage{amsfonts,amssymb}
\allowdisplaybreaks

\newcommand{\be}{\begin{equation}}
\newcommand{\ef}{\end{equation}}
\chardef\bslash=`\\ % p. 424, TeXbook
\newcommand{\G}{\Gamma}
\newcommand{\wt}{\widetilde}
\newcommand{\wh}{\widehat}
\newcommand{\ov}{\overline}
 \renewcommand{\sectionmark}[1]{}

\newcommand{\iy}{\infty}
\newcommand{\Be}{Beltrami}
\newcommand{\hol} {holomorphic}
\newcommand{\qc} {quasiconformal}

\newcommand{\fc}{\frac}

 \usepackage{amsfonts}
\newcommand{\field}[1]{\mathbb{#1}}
\newcommand{\g}{\gamma}

\newcommand{\D}{\Delta}

\newcommand{\z}{\zeta}
\newcommand{\vp}{\varphi}
\newcommand{\hC}{\widehat{\field{C}}}
\newcommand{\C}{\field{C}}
\newcommand{\R}{\field{R}}

\newcommand{\B}{\mathbf{B}}
\newcommand{\T}{\mathbf{T}}

\newcommand{\uTs} {universal Teichm\"{u}ller space}
\newcommand{\Belt} {\operatorname{Belt}}

\begin{document}

\title{Holomorphic contractibility of Teichm\"{u}ller spaces}

\author{Samuel L. Krushkal}

\begin{abstract} The problem of \hol \ contractibilty of
Teichm\"{u}ller spaces $\T(0, n)$ of the punctures spheres ($n > 4$)
arose in the 1970s in connection with solving the algebraic
equations in Banach algebras. We provide a  positive solution of
this problem.
\end{abstract}

\date{\today\hskip4mm ({\tt contractinterpol.tex})}

\maketitle

\bigskip

{\small {\textbf {2010 Mathematics Subject Classification:} Primary:
30C55, 30F60; Secondary: 30F35, 46G20}

\medskip

\textbf{Key words and phrases:} Teichm\"{u}ller spaces, Fuchsian
group, quasiconformal deformations, holomorphic contractibility,
univalent function, Schwarzian derivative, holomorphic sections

\bigskip

\markboth{Samuel Krushkal}{Holomorphic contractibility of
Teichm\"{u}ller spaces} \pagestyle{headings}

\bigskip\bigskip\centerline
{\bf 1. INTRODUCTORY REMARKS AND STATEMENT OF RESULTS}

\bigskip
Let $X$ be a complex Banach manifold which is contractible to its point
$x_0$, that is, there exists a continuous map
$F: X \times [0, 1] \to X$ with $F(x, 0) = x$ and $F(x, 1) = x_0$
for all $x \in X$.
If the map $F$ can be chosen so that for every $t \in [0, 1]$ the map
$F_t: x \mapsto F(x, t)$ of $X$ to itself is \hol \ and $F_t(x_0) = x_0$,
then $X$ is called {\em holomorphically contractible} to $x_0$.

The \hol \ contractibility of the finite dimensional Stein manifolds
closely relates to the classical Oka-Grauert $h$-principle.

First notrivial example of a (unbounded) contractible domain in
$\C^2$, which does not be holomorphically contractible, was provided
by Hirchkowitz \cite{Hi}.  Zaidenberg and Lin \cite{ZL1},
\cite{ZL2} (see also \cite{Za}) have established that there exist
contractible bounded domains of holomorphy in $\C^n, \ n > 1$, that
are not holomorphically contractible. All these domains are the
polynomial polyhedrons .

This fact became underlying for the problem of holomorphic
contractibility of Teichm\"{u}ller spaces  $\T(0, n)$ of the spheres
with $n > 4$ punctures (the case $n = 4$ is trivial, since $\T(0,
4)$ is conformally equivalent to unit disk).

This problem arose many years ago in connection with solving
algebraic equations in commutative Banach algebras and goes back to
Gorin (see, e.g., \cite{Go}). It still remains open for any
Teichm\"{u}ller space of dimension greater than $1$ (which are
topologically contractible).

The simplest example of holomorphically contractible domains in
complex Banach spaces is given by starlike domains. However all
Teichm\"{u}ller \ spaces of sufficiently great dimensions are not
stalike (see \cite{Kr2}, \cite{Kr3}, \cite{To}).

Earle \cite{Ea2} established the \hol \ contractibility of two
modified Teichm\"{u}ller spaces related to asymptotically conformal
maps.

We show that the solution of this problem is positive:

\bigskip\noindent
{\bf Theorem}. {\em Any space $\T(0, n)$ with $n > 4$ is
holomorphically contractible.}

\bigskip
As a simple consequence of this theorem, one obtains the \hol \
contractibility of two Teichm\"{u}ller spaces of Riemann surfaces
$X$ of positive genus, because the space  $\T(0, 5)$ is
biholomorphically equivalent to the space $\T(1, 2)$ of twice
punctured tori and $\T(0, 6)$ is equivalent to the space $\T(2, 0)$
of closed Riemann surfaces of genus $2$ (see, e.g., \cite{Be2}). We
present this fact in the following

\bigskip\noindent
{\bf Corollary}. {\em The spaces $\T(1, 2)$ and $\T(2, 0)$ are
holomorphically contractible.}

\bigskip
Note that the proof of the theorem involves certain specific
features of spaces $\T(0, n)$ and that the isomorphisms $\T(1, 2)
\simeq \T(0, 5)$ and $\T(2, 0) \simeq \T(0, 6)$ are exceptional.

\bigskip\bigskip
\centerline{\bf 2. A GLIMPSE OF TEICHM\"{U}LLER SPACES}

\bigskip
We briefly recall some needed facts from the Teichm\"{u}ller space
theory.

\bigskip\noindent
$\mathbf{2.1}$. \
Consider the ordered $n$-tuples of points
\begin{equation}\label{1}
\mathbf a = (0, 1, a_1, \dots, a_{n-3}, \iy), \quad n > 4,
\end{equation}
with distinct $a_j \in \C \setminus \{0, 1\}$ and the corresponding
punctured spheres
$$
X_{\mathbf a} = \hC \setminus \{0, 1, a_1 \dots, a_{n-3}, \iy\},
\quad \hC = \C \cup \{\iy\},
$$
regarded as the Riemann surfaces of genus zero. Fix a collection
$\mathbf a^0 = (0, 1, a_1^0, \dots, a_{n-3}^0, \iy)$ with $1 < a_1^0
< \dots < a_n^0 < \iy$ defining the base point $X_{\mathbf a^0}$ of
Teichm\"{u}ller space $\T(0, n) = \T(X_{\mathbf a^0})$. Its points
are the equivalence classes $[\mu]$ of \Be \ coefficients from the
ball
$$
\Belt(\C)_1 = \{\mu \in L_\iy(\C): \ \|\mu\|_\iy < 1\},
$$
under the relation that $\mu_1 \sim \mu_2$ if the corresponding
\qc \ homeomorphisms
$w^{\mu_1}, w^{\mu_2}: \ X_{\mathbf a^0} \to X_{\mathbf a}$
(the solutions of the \Be \ equation
$\overline{\partial} w = \mu \partial w$ with $\mu = \mu_1, \mu_2$)
are homotopic on $X_{\mathbf a^0}$ (and hence coincide
in the points $0, 1, a_1^0, \dots, a_{n-3}^0, \iy$).
This models $\T(0, n)$ as the quotient space
$$
\T(0, n) = \Belt(\C)_1/\sim
$$
with complex Banach structure of dimension $n - 3$ inherited from
the ball $\Belt(\C)_1$. Note that $\T(0, n)$ is a complete metric
space with intrinsic Teichm\"{u}ller metric defined by \qc \ maps.
By Royden's theorem, this metric equals the Kobayashi metric
determined by the complex structure.

Another canonical model of $\T(0, n) = \T(X_{\mathbf a^0})$ is
obtained using the uniformization of Riemann surfaces and the
holomorphic Bers embedding of Teichm\"{u}ller \ spaces. Consider the
upper and lower half-planes
$$
U  = \{z = x + i y: \ y > 0\}, \quad U^* = \{z \in \hC: \ y < 0\}
$$
and the ball
$$
\Belt(U)_1 = \{\mu \in L_\iy(\C): \ \mu|U^* = 0, \ \|\mu\|_\iy < 1\},
$$
and call the \Be \ coefficients $\mu_1$ and $\mu_2$ from this ball
{\em equivalent} if
$w^{\mu_1} = w^{\mu_2}$ on the real axes $\R = \partial U^*$
(hence on $\overline{U^*}$). Such equivalence classes $[\mu]$
are the points of the {\em \uTs} \ $T$ and correspond one-to-one to
the {\em Schwarzian derivatives}
$$
S_w(z) = \frac{w^{\prime\prime\prime}(z)}{w^\prime(z)} - \frac{3}{2}
\Bigl(\fc{w^{\prime\prime}(z)}{w^\prime(z)}\Bigr)^2
$$
of maps $w = w^\mu$ in $U^*$. These derivatives form a bounded domain
in the complex Banach space $\B = \B(U^*)$ of hyperbolically
bounded holomorphic functions $\vp$ in the lower half plane with norm
$$
\|\vp\|= \sup_{\D^*} 4 y^2 |\vp(z)|.
$$
This domain is contained in the ball $\{\|\vp\|_\B < 6\}$.

The map  $\phi: \ \mu \to S_{w^\mu}$ is \hol \ and descends to a
biholomorphic map of the space $\T$ onto this domain, which we will
identify with $\T$. It contains as complex submanifolds the
Teichm\"{u}ller spaces of all hyperbolic Riemann surfaces and of
Fuchsian groups.

\bigskip\noindent
$\mathbf{2.2}$. \
Using the \hol \ universal covering map $h: \ U \to X_{\mathbf a^0}$,
one represents the surface $X_{\mathbf a^0}$ as the quotient space
$U/\G_0$ (up to conformal equivalence), where $\G_0$ is a torsion
free Fuchsian
group of the first kind acting discontinuously on $U \cup U^*$.
The functions $\mu \in L_\iy(X_{\mathbf a^0})$ are lifted to $U$ as the
\Be \ $(-1, 1)$-measurable forms  $\wt \mu d\ov{z}/dz$ in $U$ with
respect to $\G_0$ which satisfy
$(\wt \mu \circ \g) \ov{\g^\prime}/\g^\prime = \wt \mu, \ \g \in
\G_0$ and form the Banach space $L_\iy(U, \G_0)$.

We extend these $\wt \mu$ by zero to $U^*$ and consider the unit
ball $\Belt(U, \G_0)_1$ of $L_\iy(U, \G_0)$. Then the corresponding
Schwarzians $S_{w^{\wt \mu}|U^*}$ belong to the \uTs \ $\T$ and the
subspace of such Schwarzians is regarded as the {\em Teichm\"{u}ller
space $\T(\G_0)$ of the group $\G_0$}. It is canonically isomorphic
to the space $\T(X_{\mathbf a^0})$. Moreover,
\begin{equation}\label{2}
\T(\G_0) = \T \cap \B(\G_0),
\end{equation}
where $\B(\G_0)$ is a $(n - 3)$-dimensional subspace of $\B$ which
consists of elements $\vp \in \B$ satisfying
$$
(\vp \circ \g) (\g^\prime)^2 = \vp \ \ \text{for all} \ \ \g \in \G_0
$$
(\hol \ $\G_0$-automorphic forms of degree $- 4$); see, e.g. \cite{Le}.
This leads to the representation of the space $\T(X_{\mathbf a^0})$ as a
bounded domain in the complex Euclidean space $\C^{n-3}$.

Note also that the space $\B$ is dual to the subspace $A_1(U^*)$ in
$L_1(U^*)$ formed by integrable \hol \ functions in $U^*$, while
$B(\G_0)$ has the same elements as the space $A_1(U^*, \G_0)$ of
integrable \hol \ forms of degree $- 4$ with norm $\|\vp\| =
\iint_{U^*/\G_0} |\vp(z)| dx dy$.

\bigskip\bigskip
\centerline {\bf 3. PROOF OF THEOREM}

\bigskip\noindent
$\mathbf 1^0$. \ We precede the proof of the theorem with several
lemmas which follow \cite{Kr4}.

First observe that collections (1) fill a domain $\mathcal D_n$ in
$\C^{n-3}$ obtained by deleting from this space the hyperplanes $\{z
= (z_1, \dots, z_{n-3}): \ z_j = z_l, \ j \ne l$, and with $z_1 = 0,
z_2 = 1$. This domain represents the Torelli space of the spheres
$X_{\mathbf a}$ and is covered by $\T(0, n)$. Namely, we have (cf.
e.g., \cite{Ka}; \cite{Na}, Section 2.8)

\bigskip\noindent
{\bf Lemma 1}. {\em The \hol \ universal covering space of $\mathcal
D_n$ is the Teichm\"{u}ller \ space $\T(0, n)$. This means that for
each punctured sphere $X_{\mathbf a}$ there is a \hol \ universal
covering
$$
\pi_{\mathbf a}: \T(0, n) = \T(X_{\mathbf a}) \to \mathcal D_n.
$$
The covering map $\pi_a$ is well defined by
$$
\pi_{\mathbf a} \circ \phi_{\mathbf a}(\mu) =
(0, 1, w^\mu(a_1), \dots, w^\mu(a_{n-3}), \iy),
$$
where $\phi_{\mathbf a}$ denotes the canonical projection of the ball
$\Belt(U)_1$ onto the space $\T(X_{\mathbf a})$. }

\bigskip
This lemma yields also that the truncated collections $\mathbf a_{*}
= (a_1, \dots, a_{n-3})$ provide the local complex coordinates on
the space $\T(0, n)$ and define its complex structure.

These coordinates are simply connected with the Bers local complex
coordinates on $\T(0, n)$ (related to basises of the tangent spaces
to $\T(0, n)$ at its points, see \cite{Be1}) via standard variation
of \qc \ maps of $X_{\mathbf a} = U/\G_{\mathbf a}$ (see, e.g.,
\cite{Kr1})
$$
\begin{aligned}
w^\mu(z) &= z - \fc{z(z-1)}{\pi}\iint\limits_\C
\fc{\mu(\zeta)}{\zeta(\zeta-1)(\zeta-z)} d\xi d\eta +
O(\|\mu\|_\iy^2)  \\
&= z - \fc{z(z - 1)}{\pi} \sum\limits_{\g \in \G_{\mathbf a}} \
\iint\limits_{U/\G_{\mathbf a}} \fc{\mu(\g
\zeta)|\g^\prime(\zeta)|^2}{\g \zeta(\g \zeta-1)(\g \zeta-z)} d\xi
d\eta + O(\|\mu\|_\iy^2).
\end{aligned}
$$

Now consider the ball $\Belt(U)_1$ and call its elements $\mu$
defining the same point of the \uTs \ {\em $\T$-equivalent}. The
corresponding homeomorphisms $w^\mu$ coincide on the unit circle.

We now define on this ball another equivalence relation, letting
$\mu, \ \nu \in \Belt(U)_1$ be equivalent if $w^\mu(a_j^0)
=w^\nu(a_j^0)$ for all $j$ and the homeomorphisms $w^\mu, \ w^\nu$
are homotopic on the punctured sphere $X_{\mathbf a^0}$. Let us call
such $\mu$ and $\nu$ {\em strongly $n$-equivalent}.

\bigskip\noindent
{\bf Lemma 2}. {\em If the coefficients $\mu, \nu \in \Belt(U)_1$
are $\mathbf T$-equivalent, then they are also strongly
$n$-equivalent. }

The proof of this lemma is given in \cite{Ga}.

In view of Lemmas 1 and 2, the above factorizations of the ball
$\Belt(U)_1$ generate (by descending to the equivalence classes)
a \hol \ map $\chi$ of the underlying space $\T$ into
$\T(0, n) = \T(X_{\mathbf{a^0}})$.

This map is a split immersion, i.e., it has local \hol \ sections.
In fact, we have much more:

\bigskip\noindent
{\bf Lemma 3}. {\em The map $\chi$ is surjective and has a global
\hol \ section  $s: \T(X_{\mathbf{a^0}}) \to \T$.}

\bigskip\noindent
\textbf{Proof}. The surjectivity of $\chi$ is a consequence of the
following interpolation result from \cite{CHMG}.

\bigskip\noindent
{\bf Lemma 4}. {\em Given two cyclically ordered collections of
points
$(z_1, \dots, z_m)$ and $(\z_1, \dots, \z_m)$ on the unit circle
$S^1 = \{|z| = 1\}$, there exists a \hol \ univalent function $f$
in the closure of the unit disk $\D = \{|z| < 1\}$ such that
$|f(z)| < 1$ for $z \in \ov \D$ distinct from $z_1, \dots, z_m$,
and $f(z_k) = \z_k$ for all $k = 1, \dots, m$. Moreover, there exist
univalent polynomials $f$ with such an interpolation property. }

\bigskip
Since the interpolating function $f$ given by this lemma is regular
up to the boundary, it can be extended quasiconformally across the
boundary circle $S^1$ to the whole sphere $\hC$. Hence, given a
cyclically ordered collection $(z_1, \dots, z_m)$ of points on
$S^1$, then for any ordered collection $(\z_1, \dots, \z_m)$ in
$\hC$, there is a \qc \ homeomorphism $\wh f$ of the whole sphere
$\hC$ carrying the points $z_j$ to $\z_j, \ j = 1, \dots, m$, and
such that its restriction to the closed disk $\ov \D$ is
biholomorphic on $\ov \D$.

Applying Lemma 1, one constructs \qc \ extensions of $f$ lying in
prescribed homotopy classes of homeomorphisms $X_{\mathbf z} \to
X_{\mathbf w}$. The case of maps conformal in $U$ follows from above
by conjugating the interpolating functions $f$ by the M\"{o}bius
transformation $\z \mapsto i (1 + \z)/(1 - \z)$ mapping the disk
$\Delta$ onto the lower half-plane.

\bigskip
To prove the assertion of Lemma 3 on holomorphic section for $\chi$,
take a dense subset
$$
e = \{x_1, \ x_2, \ \dots\} \subset X_{\mathbf a^0} \cap \R
$$
accumulating to all points of $\R$ and consider the surfaces
$$
X_{\mathbf a^0}^m = X_{\mathbf a^0} \setminus \{x_1, \dots, x_m\},
\quad m \ge 1
$$
(having type $(0, n + m)$). The equivalence relations on
$\Belt(\C)_1$ for $X_{\mathbf a^0}^m$ and $X_{\mathbf a^0}$ generate
a holomorphic map
$$
\chi_m: \T(X_{\mathbf a^0}^m) \to \T(X_{\mathbf a^0}).
$$

Indeed, similar to Lemma 2, we have: {\it if the coefficients $\mu,
\nu \in \Belt(U)_1$ are strongly $(n + m)$-equivalent (i.e.,
homotopic on $X_{\mathbf a^0}^m$), then they are also strongly
$n$-equivalent (homotopic on $X_{\mathbf a^0}$).}

The needed homotopy on $X_{\mathbf a^0}$ is constructed in a
standard way, for example, using the Ahlfors homotopy, letting $f(z,
t)$ be the projection of the point on the noneuclidian segment
between the corresponding covers of $f^\mu(z)$ and $f^\nu(z)$ on
hyperbolic plane $U$ which divides this segment in the proportion
$t:(1 -t)$; this homotopy extends to omitting punctures $x_j$,
together with $f^\mu$ and $f^\nu$ (cf. \cite{Ah}, \cite{Be2}).

The inclusion map $j_m: \ X_{\mathbf a^0}^m \hookrightarrow
X_{\mathbf a^0}$ forgetting the additional punctures generates a
holomorphic embedding $s_m: \ \T(X_{\mathbf a^0}) \hookrightarrow
\T(X_{\mathbf a^0}^m)$ inverting $\chi_m$.

To present this section analytically, we uniformize the surface
$X_{\mathbf a^0}^m$ by a torsion free Fuchsian group $\Gamma_0^m$ on
$U \cup U^*$ so that $X_{\mathbf a^0}^m = U/\Gamma_0^m$. By (2), its
Teichm\"{u}ller space $\T(\Gamma_0^m) = \T \cap \B(\Gamma_0^m)$. It
also can be regarded as a \hol \ universal cover of $\mathcal
D_{n+m}$.

The holomorphic universal covering maps $h: \ U^* \to U^*/\Gamma_0$
and $h^m: \ U^* \to U^*/\Gamma_0^m$ are related by $j \circ h^m = h
\circ \wh j$, where $\wh j$ is the lift of $j$. This induces a
surjective homomorphism of the covering groups $\theta_m: \Gamma_0^m
\to \Gamma_0$ by
$$
 \wh j \circ \g = \theta_m(\g) \circ \g, \quad \g \in \Gamma_0^m,
$$
and the norm preserving isomorphism $\wh j_{m,{*}}: \ \B(\Gamma_0)
\to \B(\Gamma_0^m)$ by
\begin{equation}\label{3}
\wh j_{m,{*}} \vp = (\vp \circ \wh j) (\wh j^\prime)^2,
\end{equation}
which projects to the surfaces $X_{\mathbf a^0}$ and $X_{\mathbf
a^0}^m$ as the inclusion of the space $Q(X_{\mathbf a^0})$ of
quadratic differentials corresponding to $\B(\Gamma_0)$ into the
space $Q(X_{\mathbf a^0}^m)$ (cf. \cite{EK}). The equality (3)
represents the section $s_m$ indicated above.

\bigskip To investigate the limit function for $m \to \infty$, we
embed $\T$ into the space $\B$ and compose each $s_m$ with a
biholomorphism
$$
\eta_m: \ \T(X_{\mathbf a^0}^m) \to \T(\Gamma_0^m) = \T \cap
\B(\Gamma_0^m) \quad (m = 1, 2, \dots).
$$
Then the elements of $\T(\Gamma_0^m)$ are represented in the form
$$
\wh s_m(z, \cdot) = S_{f^m}(z; X_{\mathbf a}),
$$
being parameterized by the points of $\T(X_{\mathbf a^0})$.

Each $\Gamma_0^m$ is the covering group of the universal cover $h_m:
\ U^* \to X_{\mathbf a_0^m}$, which can be normalized (conjugating
appropriately $\Gamma_0^m$) by $h_m(- i) = - i, \ h_m^\prime(- i) >
0$. Take its fundamental polygon $P_m$ obtained as the union of the
circular $m$-gon in $\Delta^*$ centered at the infinite point with
the zero angles at the vertices and its reflection with respect to
one of the boundary arcs. These polygons increasingly exhaust the
half-plane $U^*$ from inside; hence, by the Carath\'{e}odory kernel
theorem, the maps $h_m$ converge to the identity map locally
uniformly in $U^*$.

Since the set of punctures $e$ is dense on $\R$, it completely
determines the equivalence classes $[w^\mu]$ and $S_{w^\mu}$ of
$\T$, and the limit function $s(z, \cdot) = \lim_{m\to \iy} \wh
s_m(z, \cdot)$ maps $\T(X_{\mathbf a^0})$ into $\T$, what
canonically distinguishes a representative in each inverse image
$\chi^{-1}(X_{\mathbf a}) \subset \T$.

For any fixed $X_{\mathbf a}$, this function is holomorphic on
$U^*$; hence, by the well-known property of elements in the
functional spaces with sup-norms, $s(z, \cdot)$ is holomorphic also
in the norm of $\B$. This $s$ determines a holomorphic section of
the original map $\chi$, which completes the proof of Lemma 3.

\bigskip The holomorphy property indicated above is based on the
following lemma of Earle \cite{Ea1}.

\bigskip\noindent
{\bf Lemma 5}. {\em Let $E, T$ be open subsets of complex Banach
spaces $X, Y$ and $B(E)$ be a Banach space of holomorphic functions
on $E$ with sup norm. If $\vp(x, t)$ is a bounded map $E \times T
\to B(E)$ such that $t \mapsto \vp(x, t)$ is holomorphic for each $x
\in E$, then the map $\vp$ is holomorphic.}

\bigskip Holomorphy of $\vp(x, t)$ in $t$ for fixed $x$ implies the
existence of complex directional derivatives
$$
\vp_t^\prime(x,t) = \lim\limits_{\z\to 0} \fc{\vp(x, t + \z v) -
\vp(x, t)}{\z} = \fc{1}{2 \pi i} \int\limits_{|\xi|=1} \fc{\vp(x, t
+ \xi v)}{\xi^2} d \xi,
$$
while the boundedness of $\vp$ in sup norm provides the uniform
estimate
$$
\|\vp(x, t + c \z v) - \vp(x, t) - \vp_t^\prime(x,t) c v\|_{B(E)}
\le M |c|^2,
$$
for sufficiently small $|c|$ and $\|v\|_Y$.

\bigskip
The image $s(\T(X_{\mathbf a^0}))$ is an $(n-3)$-dimensional complex
submanifold in $\T$ biholomorphically equivalent to $\T(\G_0)$.

\bigskip\noindent
$\mathbf 2^0$. \ We may now prove the theorem. Pick a collection
$\mathbf a^0 = (0, 1, a_1^0, \dots, a_{n-3}^0, \iy)$ and the marked
surface $X_{\mathbf a^0}$ as indicated above, and consider its
Teichm\"{u}ller spaces $\T(X_{\mathbf a^0})$ and $\T(\G_0)$.

We embed the space $\T(0, n) = \T(X_{\mathbf a^0})$ via $\T(\G_0)$
in $\T$ and define on the space $\T(\G_0)$ a holomorphic homotopy
using the maps
$$
W^\mu = \sigma^{-1} \circ w^\mu \circ \sigma, \quad \mu \in \Belt(U)_1;
\quad \sigma(\z) = i(1 + \z)/(1 - \z), \ \ \z \in \D,
$$
and
$$
W^\mu_t(\z) := W^\mu(\z, t) = \frac{W^\mu(t \z)}{W^\mu(t)}: \ \D
\times \overline{\D} \to \hC;
$$
then
$$
w_t^\mu(z) := w^\mu(z, t) = \sigma \circ W_t^\mu \circ
\sigma^{-1}(z).
$$
By the chain rule for the Schwarzians,
 \begin{equation}\label{4}
S_{w^\mu}(\cdot, t) = t^2 S_{w^\mu}(\cdot)
=  t^{-2} (S_{W^\mu} \circ \ \sigma^{-1}) (\sigma^\prime)^{-2}.
\end{equation}
This point-wise equality determines a family of maps $\eta(\vp, t) =
S_{w^\mu_t}$ of the space $\T$ into itself, parametrized by $t \in
\overline{\D}$, with
$$
\eta(\mathbf 0, t) = \mathbf 0, \ \ \eta(\vp, 0) = \mathbf 0,
\ \ \eta(\vp, 1) = \vp.
$$
For any fixed $t$ with $|t| < 1$, the function $\eta(\vp, t)$ is
holomorphic in $\vp$ on $\T$ and by Lemma 5 for any fixed $\vp$ it
is holomorphic in $t$ in the disk $\{|t| < 1\}$ . In addition, this
function is bounded on $\T$ which follows from the estimate
$$
S_{W_t^\mu}(\z) < 6|t|^2/(|\z|^2 - 1)^2, \quad \z \in U^*.
$$
Hence, by Hartogs' theorem extended to complex Banach spaces, the
function  $\eta(\vp, t)$ is jointly holomorphic in both variables
the function $(\vp, t) \in \T \times \D$.

We apply the homotopy (4) to $\vp = S_{w^\mu} \in \T(\G_0)$. Since
it is not compatible with the group $\G_0$, there are images $\vp_t
:= \eta(\vp, t) = S_{w_t^\mu}$ which are located in $\T$ outside of
$\T(\G_0)$. The map $\chi \circ \eta(\vp, t)$ carries these images
to the points of the space $\T(0, n) = \T(X_{\mathbf a^0})$. We
compose this map with the section $s$ given by Lemma 3 and with a
biholomorphism $\xi: \ s(\T(X_{\mathbf a^0})) \to \T(\G_0)$, getting
the function
 \begin{equation}\label{5}
\Theta(\vp, t) = \xi \circ s \circ \chi \circ \eta(\vp, t)
\end{equation}
which maps holomorphically $\T(\G_0) \times \D$ into $\T(\G_0)$ with
$\Theta(\vp, 0) = \mathbf 0$.

The crucial point of the proof is to establish that the function (5)
extends holomorphically to the limit points $(\vp, 1)$ representing
the initial Schwarzians $S_{w^\mu}$. This property does not extend
(in $\B$-norm) to all points of $\T$.

To prove the limit holomorphy, fix a point $\vp_0 \in \T(\G_0)$ and
consider in its small neighborhood $V_0$ the local coordinates
$\mathbf a_{*} = (a_1, \dots, a_{n-3})$ introduced above.

Both maps $\eta$ and $\Theta$ are holomorphic in the points $(\vp_0,
t)$ of this neighborhood for all $t$ with $|t| < 1$. On the other
hand, the coordinates $\mathbf a_{*}$ are determined by the
corresponding \qc \ maps $w_t^\mu$ and, together with these maps,
are uniformly continuous in $t$ in the closed disk $\{|t| \le 1\}$.
This follows from the uniform boundedness of dilatations given by
the estimate
 \begin{equation}\label{6}
k(w_t^\mu) = \|\mu_t\|_\iy \le |t| \|\mu\|_\iy < 1
\end{equation}
(which holds for generic holomorphic motions) and from
non-increasing the Kobayashi metric $d_X(\cdot, \cdot)$ under
holomorphic maps. Since this metric on Teichm\"{u}ller spaces equals
their intrinsic Teichm\"{u}ller metric $\tau_{\T(\G_0)}$, one gets
from (6),
$$
\tau_{\T(\G_0)}(\mathbf 0, \Theta(\vp, t)) = d_{\T(\G_0)}(\mathbf 0,
\Theta(\vp, t)) \le \tanh^{-1} (|t| \|\mu\|_\iy).
$$
Hence, the function $\Theta(\vp, t)$ determines a normal family on
$V_0 \cap \T(\G_0)$.

Applying the classical Weierstrass theorem about the locally
uniformly convergent sequences of \hol \ functions in finite
dimensional domains, one derives that the limit function
$$
\Theta(\vp, 1) = \lim\limits_{t\to 1} \Theta(\vp, t)
$$
also is \hol \ on $V_0 \cap \T(\G_0)$, and then on $\T(\G_0)$, which
completes the proof of the theorem.

\medskip
{\small\em{
\leftline{Department of Mathematics, Bar-Ilan
University,}
\leftline{5290002 Ramat-Gan, Israel}
\leftline{and Department of Mathematics, University of Virginia,}
\leftline{Charlottesville, VA 22904-4137, USA}}}

\end{document}